\documentclass[amsart]{article}

\usepackage{amsmath}
\usepackage{amssymb}
\usepackage{amsfonts}
\usepackage{amsopn}
\usepackage{amsthm}
\usepackage[all]{xy}
\usepackage{hyperref} 

\swapnumbers
\theoremstyle{plain}

\newtheorem{thm}{Theorem}[section]
\newtheorem{lem}[thm]{Lemma}
\newtheorem{cor}[thm]{Corollary}

\newtheorem{prop}[thm]{Proposition}

\theoremstyle{definition}
\newtheorem{ntt}[thm]{}
\newtheorem{ex}[thm]{Example}

\newtheorem{dfn}[thm]{Definition}

\newcommand{\PP}{\mathbb{P}}  
\newcommand{\ZZ}{\mathbb{Z}}  
\newcommand{\LL}{\mathbb{L}}  
\newcommand{\OO}{\mathcal{O}} 

\newcommand{\xra}[1]{\xrightarrow{#1}}  

\newcommand{\cd}{\mathrm{cdim}}  

\newcommand{\td}{\mathrm{Td}}   

\newcommand{\CH}{\mathrm{CH}}   
\newcommand{\KK}{\mathcal{K}}   
\newcommand{\OM}{\Omega}        

\title{Degree formula for connective $K$-theory}

\author{K.~Zainoulline}

\begin{document}

\maketitle

\begin{abstract}
We apply the degree formula for connective $K$-theory
to study rational contractions of algebraic varieties.
Examples include
rationally connected varieties and complete intersections.
\footnote{E-mail: kirill@mathematik.uni-muenchen.de
}
\end{abstract}


\section{Introduction}

The celebrated {\em Rost degree formula} says that 
given a rational map $f\colon Y\dasharrow X$ 
between two smooth projective varieties 
there is the congruence relation (see \cite{Me01})
\begin{equation}\label{degfor}
\eta_p(Y)\equiv \deg f\cdot \eta_p(X)\mod n_X,
\end{equation}
where $p$ is a prime, 
$\eta_p(X)$ is the Rost number of $X$, 
$\deg f$ is the degree of $f$ and 
$n_X$ is the greatest common divisor 
of degrees of all closed points on $X$ 
(see \cite{Me03}, \cite{Mel}, \cite{Me02} and \cite{Me01}).

It was conjectured by Rost that 
the degree formula \eqref{degfor} 
should follow from a {\em generalized degree formula}
for some universal cohomology theory.
This conjecture was proven by Levine and Morel in \cite{Levine07}, 
where they constructed the theory of {\em algebraic cobordism} $\OM$
and provided the respective degree formula 
(see \cite[Theorem~1.2.14]{Levine07}).

Unfortunately, 
the generalized degree formula has one disadvantage: 
it deals with elements in the cobordism ring
which is too big and usually is hard to compute.
On the other hand the classical degree formula 
\eqref{degfor} is easy to apply
but it catches only ``pro-$p$'' effects.
The reasonable question would be to find a cohomology theory 
together with a degree formula 
which doesn't loose much information and 
is still computable.

The natural candidate for such a theory is 
the {\em connective $K$-theory} denoted by $\KK$.
It has two important properties: 
First, $\KK$ is the universal oriented cohomology theory 
for the Chow group $\CH$ and Grothendieck's $K^0$, 
meaning the following diagram of natural transformations
$$
\xymatrix{
 & \OM\ar[d]^{pr_{\KK}} & \\
\CH & \KK \ar[r]\ar[l] & K^0[\beta,\beta^{-1}],
} 
$$
where $\beta$ denotes the Bott element.
Second, it is the {\em universal birational theory} 
in the sense that
it preserves the fundamental classes for birational maps, 
i.e. for any proper birational $f\colon Y\to X$ 
we have $f_*(1_Y)=1_X$.

The respective degree formula for $\KK$
was predicted by Rost and Merkurjev
(see \cite[Example~11.4]{Me01}). 
In the present notes we deduce this formula 
from the generalized degree formula of Levine and Morel. 
Namely, we prove the following

\begin{lem}\label{mainlem} 
Let $f\colon Y\dasharrow X$ be a rational map 
between two smooth projective
irreducible varieties of the same dimension
over a field $k$ of characteristic $0$.
Then there exists a finite family of 
smooth projective varieties $\{Z_i\}_i$ over $k$ such that 
each $Z_i$ admits a projective birational map to 
a proper closed subvariety of $X$ and
\begin{equation}\label{maindeg}
\chi(\OO_Y)=\deg f\cdot \chi(\OO_X) + 
\sum_{i} n_i \cdot \chi(\OO_{Z_i}),
\end{equation}
where $n_i\in\ZZ$ and 
$\chi(\OO_X)$ is the Euler characteristic 
of the structure sheaf of $X$.
\end{lem}

Recall that an algebraic variety $X$ 
is called {\em incompressible} 
if any rational map $X\dasharrow X$ is dominant, 
i.e. has a dense image.
The notion of incompressibility appears to be very important 
in the study of the splitting properties of $G$-torsors, 
where $G$ is a linear algebraic group,
in computations of the essential and 
the canonical dimension of $G$ 
(see \cite{BR05}, \cite{CKM} and \cite{KM06}).
For instance, using the Rost degree formula \eqref{degfor} 
Merkurjev provided a uniform and shortend proof of 
the incompressibility of certain Severi-Brauer varieties,
involution varieties and quadrics 
(see \cite[\S5 and \S7]{Me01}).
Therefore, it is natural to expect that the formula \eqref{maindeg} 
can provide more examples of incompressible varieties.
As a demonstration of this philosophy we reduce 
the formula~\eqref{maindeg} to the following

\begin{thm}\label{maincong}
Let $X$ and $Y$ be smooth projective irreducible
varieties over a field
of characteristic $0$.
Let $n_X$ denote the greatest common divisor of degrees
of all closed points on $X$
and let $\tau_m$ denote the $m$-th denominator of the Todd genus. 
Assume there is a rational map $f\colon Y\dasharrow X$. 
Then we have the following congruence relation
\begin{equation}\label{DEGK}
\chi(\OO_Y)\cdot\tau_{\dim Y -1} 
\equiv \deg f\cdot \chi(\OO_X)\cdot\tau_{\dim Y-1}\mod n_X.
\end{equation}
\end{thm}

and as an immediate consequence we obtain

\begin{cor}\label{maincor}Let $X$ and $Y$ be as above.
If the image of the rational map $f\colon Y\dasharrow X$
is of dimension strictly less than the dimension of $Y$, then 
$n_X$ divides $\chi(\OO_Y)\cdot \tau_{\dim Y-1}$.
In particular,
\begin{itemize} 
\item[i)] $n_X$ always divides $\chi(\OO_X)\cdot\tau_{\dim X}$;
\item[ii)]
if $n_X$ doesn't divide $\chi(\OO_X)\cdot \tau_{\dim X-1}$, 
then $X$ is incompressible.
\end{itemize}
\end{cor}

The notion of a {\em rationally connected} variety $X$
has been introduced and extensively studied
by Campana, Koll\'ar, Miyaoka and Mori
(see \cite[IV.3.2]{Ko99} for the precise definition).
It can be easily seen (see Example~\ref{geomrat}) 
that for any rationally connected variety $X$ we have $\chi(\OO_X)=1$.
Then by Corollary~\ref{maincor} $n_X\mid \tau_{\dim X}$, and
if $n_X\nmid \tau_{\dim X-1}$ then
$X$ is {\em incompressible}.

The paper is organized as follow. First, we provide some preliminaries
on algebraic cobordism and connective $K$-theory. Then we prove
Lemma~\ref{mainlem} and Theorem~\ref{maincong} 
and discuss its applications~(Corollary~\ref{maincor}) 
to the question of incompressibility of algebraic varieties. 
In the last section
we relate the degree formula \eqref{DEGK} with the classical Rost
degree formula \eqref{degfor}.

\paragraph{Notation and conventions}

By $k$ we denote a field of characteristic $0$.
A variety will be a reduced and separated scheme of finite type over $k$.
By $pt$ we denote $\mathrm{Spec}\; k$.
Given a cycle $\alpha\in\CH(X)$ by $\deg \alpha$ we denote the push-forward
$p_*(\alpha)\in \ZZ$, where $p\colon X\to pt$ is the structure map.


\section{Algebraic cobordism and connective K-theory}

\begin{ntt} 
In \cite{Levine07} M.~Levine and F.~Morel 
introduced the theory of {\em algebraic cobordism} $\OM$ 
that is a contravariant functor from
the category of smooth quasi-projective varieties 
over a field $k$ of characteristic $0$ 
to the category of graded commutative rings. 
An element of codimension $i$ in $\OM(X)$ is 
the class $[f\colon Y\to X]$ of a projective map 
of pure codimension $i$ 
between smooth quasi-projective varieties $X$ and $Y$.
Given $f$ we denote by $f_*\colon \OM_i(Y)\to \OM_i(X)$ 
the induced push-forward 
and by $f^*\colon \OM^i(X)\to \OM^i(Y)$ the induced pull-back.
\end{ntt}

\begin{ntt} 
Let $h$ be an {\em oriented cohomology theory} as defined in \cite{Levine07}.
Roughly speaking, 
$h$ is a cohomological functor endowed with characteristic classes $c^h$. 
The main result of \cite{Levine07} says that $\OM$
is a universal oriented cohomology theory, 
i.e. any such $h$ admits a natural transformation of functors 
$pr_h\colon \OM\to h$ preserving the characteristic classes.

To any oriented cohomology theory $h$ one assigns
a one-dimensional commutative formal group law $F_h$ 
over the coefficient ring $h(pt)$ via
$$
c_1^h(L_1\otimes L_2)=F_h(c_1^h(L_1),c_2^h(L_2)),
$$
where $L_1$ and $L_2$ are lines bundles on $X$ and 
$c_1^h$ is the first Chern class.
For $\OM$ the respective formal group law turns to be
a {\em universal formal group law}
$$
F(x,y)=x+y+\sum_{i,j\ge 1} a_{ij}x^iy^j.
$$ 
The ring of coefficients $\OM(pt)$ is generated 
by the coefficients $a_{ij}$ of $F$ and 
coincides with the {\em Lazard ring} $\LL$.
\end{ntt}

\begin{ntt}
We recall several auxilliary facts about 
Chow groups $\CH$ and graded $K_0$:

Consider the augmentation map 
$\varepsilon_a\colon \LL\to\ZZ$ defined
by $a_{ij}\mapsto 0$.
Define a cohomology theory $\OM_a$
by $\OM_a(X)=\OM(X)\otimes_{\varepsilon_a} \ZZ$.
According to \cite{Levine07} 
\begin{itemize}
\item 
$\OM_a(X)$ coincides with the Chow group $\CH(X)$ of 
algebraic cycles on $X$ modulo rational equivalence;
\item
the natural transformation 
$pr_a\colon\OM(X)\to \OM_a(X)$ is surjective and its
kernel is generated by elements of positive dimensions $\LL_{>0}$
of the Lazard ring;
\item
$\OM_a$ is a universal theory 
for the additive formal group law $F_a(x,y)=x+y$.
\end{itemize}

Consider the map $\varepsilon_m\colon \LL\to \ZZ[\beta,\beta^{-1}]$
defined by $a_{11}\mapsto -\beta$ and 
$a_{ij}\mapsto 0$ for $(i,j)\neq (1,1)$.
Define a cohomology theory $\OM_m$ by
$\OM_m(X)=\OM(X)\otimes_{\varepsilon_m} \ZZ[\beta,\beta^{-1}]$. 
According to \cite{Levine07}
\begin{itemize}
\item $\OM_m(X)$ coincides with
$K^0(X)[\beta,\beta^{-1}]$, where $K^0(X)$ is Grothendieck's $K^0$ of $X$; 
\item
$\OM_m$ is a universal theory
for the multiplicative periodic formal group law $F_m(x,y)=x+y-\beta xy$.
\end{itemize}
\end{ntt}

\begin{ntt} 
The cohomology theory $\KK$ which will be the central object of our discussion
is a universal theory for both additive and multiplicative
periodic formal group laws.  It is called the {\em connective $K$-theory}
and is defined as $\KK(X)=\OM(X)\otimes_\varepsilon \ZZ[v]$,
where $\varepsilon\colon \LL\to \ZZ[v]$ 
is given by $a_{11}\mapsto -v$
and $a_{ij}\mapsto 0$ for $(i,j)\neq (1,1)$.
It has the following properties (see \cite[\S 4.3.3]{Levine07}): 
\begin{itemize}
\item 
The natural transformations
$\KK\to \OM_a$ and $\KK\to \OM_m$ are given by the evaluations
$v\mapsto 0$ and $v\mapsto \beta$ respectively.
Roughly speaking, $\KK$ can be viewed as a homotopy deformation between 
$\OM_a$ and $\OM_m$.
\item The natural transformations $\OM\to\KK$ and $\KK\to \OM_a=\CH$ 
are surjective.
\item The respective formal group law $F$
is the multiplicative non-periodic formal group
law $F_\KK(x,y)=x+y-vxy$, where $v$ is non-invertible.
\end{itemize}
\end{ntt}

\begin{ntt} \label{birinv}
We will extensively use the following fact 
(see \cite[Cor.4.2.5 and 4.2.7]{Levine07}):
\end{ntt}

\noindent 
The theory $\KK$ is universal 
among all oriented theories for which the {\em birational invariance}
holds, i.e. $f_*(1_Y)=1_X$ for any birational projective map 
$f\colon Y\to X$.
As a consequence, the kernel of the map $\varepsilon\colon\LL\to \ZZ[v]$
is the ideal generated by elements $[W]-[W']$, where
$W$ and $W'$ are birationally equivalent.


\section{Degree formula and Todd genus}

In the present section we prove Lemma~\ref{mainlem} and 
Theorem~\ref{maincong} of the introduction.
By $X$ and $Y$ we denote smooth projective irreducible varieties 
of the same dimension $d$ 
over a field $k$ of characteristic $0$.

\begin{ntt} 
Consider a rational morphism $f\colon Y\dasharrow X$.
Let $\bar{\Gamma}_f$ be the closure of its graph in $Y\times X$ and
$\bar{\Gamma}_f'\to \bar{\Gamma}_f$ be its resolution of singularities.
By the {\em generalized degree formula}  for the composite
$\bar{\Gamma}_f'\to Y\times X\xra{pr_2} X$ (see \cite[Thm.~1.2.14]{Levine07})
there exists a finite family of smooth projective varieties $\{Z_i\}_i$
such that each $Z_i$ admits a projective birational map $f_i$
on the proper closed subvariety of $X$
and
$$
[\bar{\Gamma}_f'\to X]=\deg f\cdot [X\xra{id} X] + 
\sum_{i} 
u_i\cdot [Z_i\xra{f_i} X],\text{ where }u_i\in\LL,
$$
where $\deg f=[k(Y):k(X)]$ if $f$ is dominant, and $\deg f=0$ otherwise.
Observe that by definition $\dim Z_i<d$ and $u_i\in\LL_{>0}$.

Applying 
the push-forward $p_*\colon \OM(X)\to \LL$ induced by the structure map
$p\colon X\to pt$ we obtain
$$
[\bar{\Gamma}_f']=\deg f\cdot [X]+ 
\sum_{i} u_i\cdot [Z_i].
$$
Then projecting on $\KK(pt)$ we obtain the following equality:
\begin{equation}\label{KDF}
[Y]_\KK= \deg f\cdot [X]_\KK + \sum_i (u_i)_\KK \cdot [Z_i]_\KK.
\end{equation}
Observe that $[Y]_\KK=[\bar{\Gamma}_f']_\KK$, since 
$Y$ and $\bar{\Gamma}_f'$ are birationally isomorphic and we have 
property~\ref{birinv}.
\end{ntt}

\paragraph{The Euler characteristic.} It turns out that the class $[X]_\KK$
can be expressed in terms of the {\em Euler characteristic} $\chi(\OO_X)$
of the structure
sheaf of $X$. Namely,

\begin{lem}\label{Euler} We have
$[X]_\KK=\chi(\OO_X)\cdot v^d$.
\end{lem}

\begin{proof} 
Consider the map $\varepsilon_m\colon \LL\to \ZZ[\beta,\beta^{-1}]$.
The image of the class $[X]$ is equal to
$[X]_m=p_*([\OO_X])\cdot \beta^d$, where $p_*$ is the push-forward
induced by the structure map $p\colon X\to pt$ and
the number $p_*([\OO_X])$ coincides with the Euler
characteristic $\chi(\OO_X)$ of the structure sheaf of $X$
(see \cite[\S10]{Me01})).
Since $\varepsilon_m$ factors through $\ZZ[v]$ and the
map $\ZZ[v]\to \ZZ[\beta,\beta^{-1}]$, $v\mapsto \beta$, 
is injective, we obtain the desired formula.
\end{proof}

\begin{cor}
If $X$ is birationally isomorphic to $Y$, 
then $\chi(\OO_X)=\chi(\OO_Y)$.
\end{cor} 

\begin{proof}
By property~\ref{birinv} of $\KK$ we have the equality 
$[X]_\KK=[Y]_\KK$.
\end{proof}

\begin{ntt} Observe that by the very definition
$$
\chi(\OO_X)=\sum_{i\ge 0} (-1)^i \dim_k H^i(X,\OO_X)
$$
and the number $\dim_k H^i(X,\OO_X)$ is known to be 
a {\em birational invariant}
for any smooth projective geometrically irreducible variety $X$
over a field $k$ of characteristic $0$. 
\end{ntt}

\begin{ntt}
Since $\dim_k H^i(X,\OO_X)$ 
doesn't depend on a base change, so is $\chi(\OO_X)$. 
Namely, 
if $X_l=X\times_k l$ is a base change by means of a field extension $l/k$, 
then $\chi(\OO_{X_l})=\chi(\OO_X)$.
\end{ntt}

\begin{ex} \label{geomrat}
Let $X$ be a {\em rationally connected} 
smooth projective variety over a field
$k$ of characteristic $0$ (see \cite[IV.3.2]{Ko99}). 
Observe that any {\em geometrically rational} or 
{\em unirational} variety (see \cite[IV.1]{Ko99})
provides an example of a rationally connected variety.

By \cite[IV.3.3 and IV.3.8]{Ko99}
we have $H^0(X,(\Omega^1_X)^{\otimes m})=0$ for every $m>0$.
By \cite[Example~15.2.14]{Fulton98} 
$
\dim_k H^0(X,(\Omega^1_X)^{\otimes m})=
\dim_k H^m(X,\OO_X)$ for every $m$
and, therefore, $\chi(\OO_X)=1$.
\end{ex}

Combining~\eqref{KDF} and Lemma~\ref{Euler} we prove
the Lemma~\ref{mainlem} of the introduction:

\begin{lem} Let $f\colon Y\dasharrow X$ be a rational map between two
smooth projective irreducible
varieties of the same dimension $d$
over a field of characteristic $0$.
Then there exists a finite family of smooth projective 
varieties $\{Z_i\}_i$ such that each $Z_i$ admits a projective birational
map $f_i$ on a proper closed subvariety of $X$ and
\begin{equation}\label{DF}
\chi(\OO_Y)=\deg f\cdot \chi(\OO_X) + 
\sum_{i} n_i \cdot \chi(\OO_{Z_i}),\text{ where }n_i\in\ZZ.
\end{equation}
\end{lem}

\paragraph{The Todd genus.} 
It is well-know that the Euler characteristic $\chi(\OO_X)$
is related with the {\em Todd genus} of $X$.

\begin{ntt} According to \cite[Example~3.2.4]{Fulton98}
the {\em Todd class} of the tangent bundle of
a smooth projective variety $X$ is the following polynomial
in Chern classes
\begin{equation}\label{tdgen}
\td(X)=1+\tfrac{c_1}{2}+\tfrac{c_1^2+c_2}{12}+\tfrac{c_1c_2}{24} + 
\tfrac{-c_1^4+4c_1^2c_2+3c_2^2+c_1c_3-c_4}{720}+\ldots,
\end{equation}
where $c_i\in \CH^i(X)$ 
denotes the $i$-th Chern class of the tangent bundle of $X$.
The denominators $\tau_0=1$, $\tau_1=2$, $\tau_2=12$, $\tau_3=24$, 
$\tau_4=720\ldots$ are called {\em Todd numbers}. 
We have the following explicit formula for $\tau_d$ 
(see \cite[Example~9.9]{Me01}):
\begin{equation}\label{explform}
\tau_d=\prod_{p\; prime} p^{[\tfrac{d}{p-1}]}.
\end{equation}
In particular, $\tau_{d-1}\mid \tau_d$ for any $d$.
\end{ntt}

\begin{ntt}
To compute the Euler characteristic $\chi(\OO_X)$ we 
may use the following equality (see \cite[Corollary~18.3.1]{Fulton98}):
\begin{equation}\label{tdchi}
\chi(\OO_X)=\deg \td(X),
\end{equation}
where $\deg \td(X)$ is the degree of the 
$d$-th homogeneous component of the Todd class $\td(X)$ and 
is called the Todd genus of $X$.

Observe that the Euler characteristic 
and the Todd genus are multiplicative, i.e.
$\chi(\OO_{X\times Y})=\chi(\OO_X)\cdot \chi(\OO_Y)$.
\end{ntt}

\begin{ex}\label{hypers} Let $X$ be a complete
intersection of $m$ smooth hypersurfaces of degrees
$d_1,\ldots,d_m$ in $\PP^n_k$. Then
$$
\chi(\OO_X)=
Res_{z=0}\frac{\prod_{i=1}^m(1-e^{-d_iz})}{(1-e^{-z})^{n+1}}.
$$
Indeed, by \cite[Example~15.2.12.(iii)]{Fulton98} we have
$$
\chi(\OO_X)=\deg \big(\td(T_{\PP^n})\cdot\prod_{i=1}^m(1-e^{-d_iz})\big)=
\deg \big( (\frac{z}{1-e^{-z}})^{n+1}\cdot \prod_{i=1}^m(1-e^{-d_iz}) \big),
$$
where $z=c_1(\OO_{\PP^n}(1))$. 
\end{ex}

We are ready now to prove the main result of this paper 
(Theorem~\ref{maincong})

\begin{thm}\label{propinc}
Let $X$ and $Y$ be smooth projective irreducible
varieties over a field $k$ 
of characteristic $0$.
Let $n_X$ denote the greatest common divisor of degrees
of all closed points on $X$ 
and let $\tau_m$ denote the $m$-th Todd number. 
Assume there is a rational map $f\colon Y\dasharrow X$. 
Then 
\begin{equation}\label{DFR}
\chi(\OO_Y)\cdot\tau_{\dim Y -1} 
\equiv \deg f\cdot \chi(\OO_X)\cdot\tau_{\dim Y-1}\mod n_X.
\end{equation}
\end{thm}

\begin{proof} Taking the product with a projective space
we may assume
that $X$ and $Y$ have the same dimension $d$.
By the formula \eqref{DF} we have
$$
\chi(\OO_Y)=
\deg f\cdot \chi(\OO_X)+\sum_{i} n_i\chi(\OO_{Z_i})=
\deg f\cdot \chi(\OO_X)+\sum_{j=0}^{d-1}\Big(
\sum_{i,\dim Z_i=j} n_i\chi(\OO_{Z_i})\Big), 
$$
where each $Z_i$ admits a projective morphism
$Z_i\to X$. Identifying the Euler characteristic with the Todd number
using \eqref{tdchi} we see that 
all the characteristic numbers in the numerator 
of the $(\dim Z_i)$-th homogeneous component of the polynomial \eqref{tdgen}
for $\td(Z_i)$  
are divisible by $n_X$.
Therefore, we have
$$
\sum_{j=0}^{d-1}\Big(
\sum_{i,\dim Z_i=j} n_i \cdot \frac{n_X\cdot m_i}{\tau_j}\Big)
=\sum_{j=0}^{d-1}\frac{n_X}{\tau_j}\Big(\sum_{i,\dim Z_i=j}n_im_i\Big)=
$$
$$
=\frac{n_X}{\tau_{d-1}}\sum_{j=0}^{d-1}
\frac{\tau_{d-1}}{\tau_j}\Big(\sum_{i,\dim Z_i=j}n_im_i\Big)=
\frac{n_X\cdot m}{\tau_{d-1}},
$$
where $m=\sum_{j=0}^{d-1}
\frac{\tau_{d-1}}{\tau_j}\big(\sum_{i,\dim Z_i=j}n_im_i\big)\in\ZZ$ according
to \eqref{explform}. This completes the proof of the theorem.
\end{proof}

\begin{ex}\label{cursur} 
If $X$ and $Y$ are curves, then 
the degree formula \eqref{DFR} turns into
\begin{equation}\label{curdeg}
\chi(\OO_Y)\equiv \deg f\cdot \chi(\OO_X) \mod n_X,
\end{equation}
where $n_X$ denotes the g.c.d.
of degrees of all closed points on $X$.
Observe that
$\chi(\OO_X)=1-p_g$, where $p_g$ is the geometric genus of a
geometrically irreducible curve $X$.

For surfaces $X$ and $Y$ it can be stated as
\begin{equation}\label{surfdeg}
2\chi(\OO_Y)\equiv \deg f\cdot 2\chi(\OO_X)\mod n_X.
\end{equation}
Observe that $\chi(\OO_X)=1-q+p_g$, where 
$q$ is the irregularity and $p_g$ is the geometric
genus of a geometrically irreducible surface $X$.
\end{ex}


\section{Incompressibility}

\begin{ntt}
The notion of (in-)compressibility of algebraic varieties 
appears naturally in the study of the splitting properties
of $G$-torsors and their canonical dimensions.
Recall that (see \cite[\S4]{KM06} and \cite[\S1]{CKM}) a {\em canonical
dimension} $\cd X$ of a smooth projective irreducible variety $X$ over $k$
is defined to be the minimal dimension 
of a closed irreducible subvariety $Z$
of $X$ such that $Z_{k(X)}$ has a rational point.
By definition we have $\cd X\le \dim X$.
If $\cd X=\dim X$, then $X$ is called {\em incompressible}.

To say that ``$Z_{k(X)}$ has a rational point'' 
is the same as to say that 
there is a rational dominant map $X\dasharrow Z$. Therefore, 
a variety $X$ is incompressible 
if and only if any rational map $X\dasharrow X$
is dominant.
\end{ntt}

As an immediate consequence of Theorem~\ref{propinc} we obtain
\begin{cor}\label{maincorr}
Let $X$ and $Y$ be smooth projective irreducible
varieties over a field $k$ 
of characteristic $0$.
Let $n_X$ denote the greatest common divisor of degrees
of all closed points on $X$ 
and let $\tau_m$ denote the $m$-th Todd number. 
Assume there is a rational map $f\colon Y\dasharrow X$. 

If the image of $f$
is of dimension strictly less than the dimension of $Y$, then 
$n_X$ divides $\chi(\OO_Y)\cdot \tau_{\dim Y-1}$.
In particular, i) $n_X\mid \chi(\OO_X)\cdot\tau_{\dim X}$, and 
ii) if $n_X\nmid \chi(\OO_X)\cdot \tau_{\dim X-1}$, then $X$ is incompressible.
\end{cor}

\begin{proof} To prove i) we apply Theorem~\ref{propinc} to the projection
$X\times\PP^1\to X$, and to prove ii) we apply~\ref{propinc} to
the identity $X\to X$.
\end{proof}

\begin{ex} Let $X$ be a rationally connected smooth projective variety
over a field $k$ of characteristic $0$.
Then according to Example~\ref{geomrat} we have $\chi(\OO_X)=1$. 
Therefore, $n_X\mid \tau_{\dim X}$ and
\begin{equation}\label{grim}
n_X\nmid \tau_{\dim X-1}\Longrightarrow X\text{ is incompressible.}
\end{equation}
Observe that if $X$ is a curve or a surface,
the implication \eqref{grim} can be proven
directly using the geometry 
(see \cite[\S8]{BR05} and \cite[\S2, \S3]{CKM}). 
\end{ex}

\begin{ex}(cf. \cite[Example~8.2]{Mel} and \cite[\S7.3]{Me03})\label{myex} 
Let $X$ be a complete intersection of 
$m$ smooth hypersurfaces 
of degrees $d_1,\ldots,d_m$ in $\PP^n_k$.

Assume that $\dim X=p-1$ for some prime $p$.
Let $m_p$ denote the number of degrees $d_i$ 
which are divisible by $p$.
We claim that 
\begin{equation}
p\nmid m_p \text{ and }p\nmid \tfrac{d_1d_2\ldots d_m}{n_X}
\Longrightarrow X\text{ is incompressible.}
\end{equation}

Indeed, by the formula from Example~\ref{hypers} the Euler characteristic 
$\chi(\OO_X)$ is equal to the coefficient at $z^n$ in the
expansion of
$$
d_1d_2\ldots d_m z^m
\prod_{i=1}^m\Big(\sum_{r=0}^{\dim X}\tfrac{(-d_i)^r}{(r+1)!}z^r\Big)
\Big(\sum_{r=0}^{\dim X}\tfrac{B_r}{r !}z^r\Big)^{p+m}, 
$$
where $B_r$ denotes the $r$-th Bernoulli number.
Since the denominator
of $\tfrac{B_r}{r!}$ is not divisible by $p$ for any $r<p-1$ and
is divisible by $p$ for $r=p-1$, we obtain
$$
\frac{\chi(\OO_X)\tau_{p-2}}{n_X}=
\frac{\tau_{p-2}\cdot d_1d_2\ldots d_m}{n_X}\cdot
\Big(\frac{a}{b}\
-\frac{m_p}{p!}\Big)\notin \ZZ,\text{ where }p\nmid ab.
$$
\end{ex}

\begin{ex} Let $Y$ be a smooth hypersurface of degree $p^r$, $r>0$,
in $\PP^p_k$ with $n_Y=p^r$. Assume that there is a rational map
$Y\dasharrow X$ to a smooth projective variety $X$ with $\dim X<\dim Y$. 
Then $n_X\mid p^{r-1}$.

By \cite[Prop.~6.2]{Ko99} we have $n_X\mid n_Y$. Therefore, $n_X\mid p^r$.
By Corollary~\ref{maincorr}
$n_X\mid \chi(\OO_Y)\cdot\tau_{p-2}$ and by the previous example the right hand
side is not divisible by $n_Y=p^r$.
\end{ex}

The next proposition provides a version of the index reduction formula
for varieties with different Euler characteristics

\begin{prop}\label{ired}
Let $X$ and $Y$ be smooth projective geometrically irreducible 
varieties over a field $k$ of characteristic $0$. Let $p$ be a prime.
Assume that 
\begin{itemize}
\item $\dim X$, $\dim Y<p$;
\item $n_X=n_Y=p$;
\item $p\mid \chi(\OO_X)$ and $p\nmid \chi(\OO_Y)$. 
\end{itemize}
Then $n_{X_{k(Y)}}=p$, where $k(Y)$ is the function field of $Y$.
\end{prop}

\begin{proof}
Taking the product with a projective space we may
assume that $\dim X=\dim Y=p-1$. 
Obviously, $n_{X_{k(Y)}}\mid n_X$.
Assume that $n_{X_{k(Y)}}=1$.
Then $X_{k(Y)}$ has a closed point $P$ of degree $m$ coprime to $p$.

We follow the proof of \cite[Thm.~3.3]{CM04} 
(see also  \cite[Thm.~5.1]{Co05}).
Let $K$ denote the residue field of $P$
and let $Y'$ be a smooth projective variety over $k$ such that
$K=k(Y')$. Such $Y'$ can be always obtained from $Y$ by taking
the normalization and then resolving the singularities.
The condition that $X_{k(Y)}$ has the point $P$
means that there is a rational map $f_P\colon Y'\dasharrow X$.

Assume that $f_P$ is not dominant. Then by Theorem~\ref{propinc} applied
to the maps $f_P\colon Y'\dasharrow X$ and $Y'\to Y$ we obtain that 
$$
n_X\mid \chi(\OO_{Y'})\cdot \tau_{p-2}\;\text{ and }\;
\chi(\OO_{Y'})\cdot \tau_{p-2}\equiv
m\cdot\chi(\OO_Y)\cdot\tau_{p-2}\mod n_Y.
$$ 
Therefore,
$p\mid \chi(\OO_Y)$, a contradiction.
Hence, $f_P$ has to be dominant.

If $f_P$ is dominant, then by~\eqref{DFR} we have
$$
\chi(\OO_Y)\cdot \tau_{p-2}\equiv
\deg(f_P)\cdot\chi(\OO_X)\cdot\tau_{p-2}\mod n_X.
$$ 
Since $p\mid \chi(\OO_X)$, the left hand side has to be divisible by $p$,
a contradiction.

Therefore, $n_{X_{k(Y)}}\neq 1$, and the proposition is proven.
\end{proof}

\begin{cor}(cf. \cite[Thm.~1.1]{Co05})\label{cohd}
Let $X/k$ be as above, i.e. $\dim X<p$, $n_X=p$ and 
$p\mid \chi(\OO_X)$ for some prime $p$. Then there exists a field extension
$K/k$ of cohomological dimension one such that $n_{X_K}=p$.
\end{cor}

\begin{proof} 
Take $Y$ to be a 
Severi-Brauer variety of a division algebra of degree $p$ and apply
Proposition~\ref{ired}. We obtain that $n_{X_{k(Y)}}=n_X=p$.
Repeating the arguments of \cite[\S 2]{Co05} 
for the given $X$ and $Y$ we finish the proof.
\end{proof}

\begin{ex} 
If $X$ is a curve of genus one (here $\chi(\OO_X)=0$) 
the existence of such a field $K$ was proven in \cite{Du98}.

Let $p$ be an odd prime and 
let $X$ be a smooth hypersurface of degree $p$ in $\PP^{p-1}_k$
with $n_X=p$. 
Observe that $X$ is an anisotropic {\em Calabi-Yau} variety with
$\chi(\OO_X)=0$. By Corollary~\ref{cohd} we obtain that
there exists a field $K/k$ of cohomological dimension one such that
$n_{X_K}=p$. In particular, $X$ has no zero-cycles of degree one over $K$.
\end{ex}


\section{Comparison with the classical degree formulas}

\begin{ntt}\label{defpart} We follow the notation of \cite{Me02}.
Let $\alpha=(\alpha_1,\alpha_2,\ldots,\alpha_r)$ be a partition, i.e.
a sequence of integers (possibly empty) 
$0<\alpha_1\le \alpha_2\le\ldots \le \alpha_r$, 
and let $|\alpha|=\alpha_1+\alpha_2+\ldots+\alpha_r$ denote its degree.
For any $\alpha$ we define the smallest symmetric polynomial
$$
P_\alpha(x_1,x_2,\ldots)=\sum_{(i_1,i_2,\ldots,i_r)}
x_{i_1}^{\alpha_1}x_{i_2}^{\alpha_2}\ldots x_{i_r}^{\alpha_r}=
Q_\alpha(\sigma_1,\sigma_2,\ldots)
$$
containing the monomial $x_1^{\alpha_1}x_2^{\alpha_2}\ldots x_r^{\alpha_r}$
where the $\sigma_i$ are the elementary symmetric functions.
\end{ntt}

\begin{ntt}
Let $X$ be a smooth projective variety of dimension $d$.
Let $c_i=c_i(-T_X)$ denote the $i$-th Chern class 
of the inverse of the tangent bundle. 
Let $\alpha$ be a partition of $d$.
We define the $\alpha$-characteristic number of $X$ by
$$
c_\alpha=\deg Q_{\alpha}(c_1,c_2,\ldots).
$$
Observe that $c_{(1,1,\ldots,1)}=\deg c_{d}(-T_X)$ and 
$c_{(d)}$ defines the so called additive characteristic number of $X$.
\end{ntt}

\begin{ntt}\label{despar} We fix a prime $p$. Consider a partition 
\begin{equation}
\alpha=(p-1,\ldots,p-1,p^2-1,\ldots,p^2-1,\ldots),
\end{equation}
where $p^i-1$ is repeated $r_i$ times (in \cite{Me03} it 
was denoted by the sequence
$R=(r_1,r_2,\ldots)$). 
The set of all such partitions $\alpha$ will
be denoted by $\Lambda_p$.
According to \cite[\S6]{Me03} for any $\alpha\in\Lambda_p$
the characteristic number $c_\alpha$ 
is divisible by $p$.
The integer $\tfrac{1}{p}c_{(p-1,\ldots,p-1)}$ is called
the {\em Rost number} and is denoted by $\eta_p$. 
\end{ntt}

\begin{ntt}
By \cite[Theorem~6.4]{Me03} for any prime $p$
and any partition $\alpha\in\Lambda_p$ of $d$
we have the degree formula:  
\begin{equation}\label{rostdeg}
\frac{c_\alpha(Y)}{p}\equiv \deg f\cdot \frac{c_\alpha(X)}{p}\mod n_X,
\end{equation}
where $\deg f$ is the degree of a rational map $f\colon Y\dasharrow X$
and $d=\dim Y=\dim X$.
In particular, if $n_X\nmid \tfrac{1}{p}c_\alpha(X)$, 
then $X$ is incompressible.
\end{ntt}

In the present section we discuss the relations between
the classical degree formulas~\eqref{rostdeg} 
and the degree formula~\eqref{DFR}. 
The following lemma provides an explicit formula for $\chi(\OO_X)$
in terms of characteristic numbers $c_\alpha(X)$

\begin{lem}\label{deschi} 
Let $X$ be a smooth projective variety of dimension $d$ over $k$.
Then 
$$
\chi(\OO_X)=(-1)^d\sum_{\alpha,|\alpha|=d} 
\frac{c_\alpha(X)}{(\alpha_1+1)!(\alpha_2+1)!\ldots (\alpha_r+1)!}.
$$
\end{lem}

\begin{proof} 
By the very definition $\td(X)=\td(T_X)=\prod_{i=1}^d Q(x_i)$,
where $Q(x_i)=\frac{x_i}{1-e^{-x_i}}$ and $x_i$ are the roots
of the tangent bundle $T_X$. 
Since $\td(T_X)=\td(-T_X)^{-1}=\prod_{i=1}^d Q(-x_i)^{-1}$, 
its component of degree $d$ is equal
to the coefficient at $z^d$ in the expansion of the product
$$
\prod_{i=1}^d 
(1-\tfrac{x_i}{2!}z+\tfrac{x_i^2}{3!}z^2-\tfrac{x_i^3}{4!}z^3+\ldots).
$$
Analyzing the product we see that this coefficient is, indeed, given by
$$
(-1)^d\sum_{\alpha,|\alpha|=d} 
\frac{P_\alpha(x_1,x_2,\ldots)}{(\alpha_1+1)!(\alpha_2+1)!\ldots},
$$
where $P_\alpha$ is the minimal symmetric polynomial from \ref{defpart}.
\end{proof}

\begin{lem}\label{prpart} 
Let $\alpha=(\alpha_1,\ldots,\alpha_r)$ be a partition
of $d$ and let $v_p(m)$ denote the $p$-adic valuation of an integer $m$. 
Then 
$$
v_p(\tau_{d-1})+1\ge v_p\Big(\prod_{i=1}^r(\alpha_i+1)!\Big),
$$ 
where the equality holds if and only if $\alpha\in\Lambda_p$.
\end{lem}

\begin{proof} Follows from the formulas 
$v_p(\tau_{d-1})=[\tfrac{d-1}{p-1}]$ and 
$v_p(m!)=\sum_{j=1}^\infty [\tfrac{m}{p^j}]$.
\end{proof}

\begin{dfn} Let $p$ be a prime $p$ and let $d$ be an integer. 
We define a linear combination $u_p$ of characteristic numbers $c_\alpha$,
$\alpha\in\Lambda_p$, as
$$
u_p=\sum_{\alpha\in\Lambda_p, |\alpha|=d} 
\frac{\tau_{d-1}}{\prod_{i=1}^r(\alpha_i+1)!}c_\alpha.
$$
According to Lemma~\ref{prpart} we have
$$
u_p=\sum_{\alpha\in\Lambda_p, |\alpha|=d}
n_\alpha \frac{c_\alpha}{p},
\text{ where }
n_\alpha=\frac{p\cdot\tau_{d-1}}{\prod_{i=1}^r (\alpha_i+1)!}\in\ZZ\text{ and }
p\nmid n_\alpha.
$$
\end{dfn}

\begin{prop} Let $X$ be a smooth projective variety of dimension $d$
over a field of characteristic $0$. Then
$$
n_X\nmid \chi(\OO_X)\tau_{d-1}\Longleftrightarrow 
\exists p\text{ such that } n_X\nmid u_p(X).
$$
\end{prop}

\begin{proof} 
Since $n_X\mid c_\alpha$ for any $\alpha$, 
$u'_p=\tfrac{pu_p}{n_X}\in\ZZ$. We have the following chain
of equivalences

$n_X\mid \chi(\OO_X)\tau_{d-1}\Longleftrightarrow
n_X\mid \sum_p \tfrac{n_X}{p}u_p'
\Longleftrightarrow
\sum_p \tfrac{u_p'}{p}\in \ZZ$

$\Longleftrightarrow \forall p\;\; p\mid u_p'
\Longleftrightarrow \forall p\;\; \tfrac{u_p}{n_X}\in\ZZ$.
\end{proof}

\begin{ex} Let $X$ be a smooth projective curve. In this case
we have only one non-trivial partition $\alpha=(1)\in\Lambda_2$ and
$$
\chi(\OO_X)=-\tfrac{1}{2}c_{(1)}(X)=-\eta_2(X)=-u_2(X).
$$
Therefore, for curves the degree
formula \eqref{DFR} coincides with the classical one.
\end{ex}

\begin{ex} 
For a smooth projective surface $X$ 
we have two partitions $\alpha=(1,1)$ and $(2)$, where
the first one belongs to $\Lambda_2$ 
and the second one to $\Lambda_3$.
We have
$$
\chi(\OO_X)=\tfrac{1}{4}c_{(1,1)}(X)+\tfrac{1}{6}c_{(2)}(X)
\text{ and }u_2=\tfrac{1}{2}c_{(1,1)}=\eta_2,\;
u_3=\tfrac{1}{3}c_{(2)}=\eta_3.
$$
The degree formula~\eqref{DFR}
turns into a sum of the classical
degree formulas
$$
(\eta_2+\eta_3)(Y)\equiv \deg f\cdot (\eta_2+\eta_3)(X) \mod n_X
$$
and
$$
n_X\nmid \tau_{d-1} \chi(\OO_X) \Longleftrightarrow
n_X\nmid \eta_2(X)\text{ or }n_X\nmid \eta_3(X).
$$
So from the point of view of incompressibility the degree formula
\eqref{DFR} provides the same answer as the classical degree formulas.
\end{ex}

\begin{ex} 
Let $X$ be a smooth projective $3$-fold. 
In this case we have three partitions $(1,1,1)$, $(1,2)$ and $(3)$,
where the first and the last one belong to $\Lambda_2$.
We have
$$
\chi(\OO_X)=-\tfrac{1}{8}c_{(1,1,1)}- \tfrac{1}{12}c_{(1,2)}-
\tfrac{1}{24}c_{(3)}\text{ and }
u_2= \tfrac{3}{2}c_{(1,1,1)}+ \tfrac{1}{2}c_{(3)}.
$$
Therefore,
$$
n_X\nmid \chi(\OO_X)\tau_2 \Longleftrightarrow 
n_X\nmid u_2
\stackrel{(*)}\Longrightarrow n_X\nmid \tfrac{1}{2}c_{(1,1,1)}\text{ or }
n_X\nmid \tfrac{1}{2}c_{(3)}
$$
which means that for $3$-folds the classical degree formulas \eqref{rostdeg} 
detect the incompressibility better than \eqref{DF}.
\end{ex}

\begin{ntt}
Since each characteristic number $c_\alpha$
is divisible by $n_X$, to say that $n_X\mid q\cdot c_\alpha$, 
$q\in \mathbb{Q}$,
is equivalent to say that $q\cdot C_\alpha\in\ZZ$, 
where $C_\alpha:=c_\alpha/n_X\in\ZZ$. 
Hence, the implication ($*$) can be rewritten as
$$
\tfrac{1}{2}C_{(1,1,1)}+ \tfrac{1}{2}C_{(3)}\not\in \ZZ\Longrightarrow
\tfrac{1}{2}C_{(1,1,1)}\not\in\ZZ\text{ or }\tfrac{1}{2}C_{(3)}\not\in\ZZ
$$
In particular, the implication ($*$) becomes an equivalence if and only if
a $3$-fold $X$ satisfies the following condition
\begin{equation}\label{cond3}
C_{(1,1,1)}=\tfrac{c_{(1,1,1)}(X)}{n_X}\text{ is even or }
C_{(3)}=\tfrac{c_{(3)}(X)}{n_X}\text{ is even.}
\end{equation}
\end{ntt}

\begin{ex}
Let $X$ be a complete intersection of $m$ hypersurfaces of degrees
$d_1,\ldots,d_m$ in $\PP^{m+3}$.
Using the formula from \cite[\S8]{Mel} we obtain:
$$
C_{(3)}=\frac{\partial}{\partial z}\Big(\frac{\prod_{i=1}^m(1+d_i^3z)}{(1+z)^{m+4}}\Big)_{z=0}
=(\sum_{i=1}^m d_i^3)-m-4\equiv \sigma_1+m\mod 2,
$$
where $\sigma_1$ denotes the sum of all degrees. And
$$
C_{(1,1,1)}=
\frac{1}{6}\frac{\partial^3}{\partial z^3}
\Big(\frac{\prod_{i=1}^m(1+d_iz)}{(1+z)^{m+4}}\Big)_{z=0}\equiv
{m+2\choose 3}+{m+1\choose 2}\sigma_1+m\sigma_2\mod 2,
$$
where $\sigma_2$ is the second elementary symmetric function
in $d_i$-s.
 
Hence, a complete intersection $X$ satisfies \eqref{cond3}
if and only if it satisfies one of the following conditions

\begin{tabular}{ll}
$\bullet$ $m=4k$; & $\bullet$ $m=4k+2$, $\sigma_1$ is even;\\
$\bullet$ $m=4k+1$, $\sigma_1$ or  $\sigma_2$ is odd; &
$\bullet$ $m=4k-1$, $\sigma_1$ is odd or $\sigma_2$ is even.
\end{tabular}
\end{ex}

\paragraph{Acknowledgements} I would like to thank 
Z.~Reichstein, Y.~Tschinkel and A.~Vishik for useful comments. 
I am very grateful to J.-L. Colliot-Th\'el\`ene for the encouraging
attention to my work and for numerous suggestions which helped to impove
the exposition of the paper.
This work was supported by SFB~701,
INTAS~05-1000008-8118 and DFG~GI~706/1-1.

\vspace{1cm}
\noindent Kirill Zainoulline,
Mathematisches Institut der LMU M\"unchen\\
Theresienstr. 39,
D-80333 M\"unchen, Germany

\end{document}